\documentclass[smallextended]{svjour3}%[12pt]{elsarticle}

\usepackage{amssymb}
\usepackage{amsmath}
\usepackage{graphicx}
\usepackage{verbatim}
\usepackage{amssymb}
\usepackage{graphicx}
\usepackage{rotating}
\usepackage{lscape}
\usepackage{subfig}
\usepackage{algorithm}
\usepackage[noend]{algpseudocode}

\usepackage{amsthm}
\usepackage{moresize}
\usepackage{setspace}
\usepackage{subfig}
\usepackage[numbers]{natbib}

\begin{document}
	
%\setlength{\abovedisplayskip}{3pt}
%\setlength{\belowdisplayskip}{3pt}

%\begin{frontmatter}



\title{QUBO transformation using Eigenvalue Decomposition}

%\author{Amit Verma\corref{cor1}}
%\ead{averma@missouriwestern.edu}
%\fntext[fn1]{Student}

%\author{Mark Lewis}
%\ead{mlewis14@missouriwestern.edu}
%\fntext[fn2]{Lecturer}

%\address{Craig School of Business, Missouri Western State University, Saint Joseph, MO, 64507, United States}
\author{Amit Verma         \and
	Mark Lewis %etc.
}


\institute{Amit Verma \at
	Craig School of Business, Missouri Western State University, Saint Joseph, MO, 64507, United States \\
	Tel.: +1 816-271-4357\\
	Fax: +1 816-271-4508\\
	\email{averma@missouriwestern.edu}           %  \\
	\and
	Mark Lewis \at
	\email{mlewis14@missouriwestern.edu}
}

\date{Received: date / Accepted: date}

\maketitle

\begin{abstract}
	
	Quadratic Unconstrained Binary Optimization (QUBO) is a general-purpose modeling framework for combinatorial optimization problems and is a requirement for quantum annealers. This paper utilizes the eigenvalue decomposition of the underlying $Q$ matrix to alter and improve the search process by extracting the information from dominant eigenvalues and eigenvectors to implicitly guide the search towards promising areas of the solution landscape. Computational results on benchmark datasets illustrate the efficacy of our routine demonstrating significant performance improvements on problems with dominant eigenvalues.


\end{abstract}

%\begin{keyword}
%Quadratic Unconstrained Binary Optimization, pseudo-Boolean optimization, eigenvalue  decomposition, eigenvalues, eigenvectors, quantum computer, maximum diversity problem%, local optima network
%\end{keyword}

%\end{frontmatter}


\section{Introduction}
%\cite{billionnet2004exact}
%STRICT INEQUALITY - Discuss
%$A^{t}$
Quadratic Unconstrained Binary Optimization (QUBO) is a popular modeling framework wherein many optimization problems are transformed into $max \; x^{t}Qx; x \in \{0,1\}$ where $Q$ is a $n X n$ symmetric matrix and $x$ denotes the binary decision vector of size $n$. QUBO has emerged into significance in the last decade due to advances in quantum annealers which rely on $0/1$ modeling format (or equivalent Ising spin $-1/+1$). For a detailed survey of applications using this modeling format, we refer the readers to \cite{kochenberger2014unconstrained}. %In most applications, the symmetric $Q$ matrix with coefficient $q_{ij}$ is transformed into an upper or lower diagonal matrix with $q_{ij}^{\prime} = 2 q_{ij}$. However, eigendecomposition benefits from the symmetric nature of the $Q$ matrix.

While the $Q$ matrix is often transformed into an upper diagonal matrix, eigendecomposition will retain a symmetric format because a symmetric matrix always leads to real eigenvalues and eigenvectors. According to spectral decomposition theorem (\cite{halmos1963does}), a symmetric matrix $Q$ can be decomposed into $\sum_{i=1}^{n} \lambda_i c_i c_i^{t}$ where $\lambda_i$ and $c_i$ denotes the $i'th$ eigenvalue and eigenvector respectively. Thus, $max \; x^{t}Qx$ can be transformed into $max \; x^{t} (\sum_{i=1}^{n} \lambda_i c_i c_i^{t}) x$. Since each $\lambda_i$ is a scalar, this yields $max \; \sum_{i=1}^{n} \lambda_i (x^{t} c_i) (c_i^{t} x)$ where $x$ and $c_i$ are vectors and $x^{t} c_i = c_i^{t} x$. Thus, we rewrite the final expression as $max \; \sum_{i=1}^{n} \lambda_i (x^{t} c_i)^2 $.

The time complexity of the complete eigendecomposition ranges from $O(n^2)$ to $O(n^3)$ depending on the structure of the matrix, matrix multiplication routine and employed algorithm, according to (\cite{pan1999complexity}). However, the top $k$ eigenvalues and associated eigenvectors could be computed more efficiently using techniques described in \cite{sorensen1997implicitly} and are part of numerical software packages like ARPACK (\cite{lehoucq1998arpack}). In this paper, we will focus on top $k$ eigenvalues and eigenvectors for the $Q$ matrix. As an example the complete eigendecomposition of the following $Q$ matrix is given by:
%$\resizebox{0.91\hsize}{!}{%
\[
\begin{bmatrix}
	-7 & 2 & 2 \\
	2 & 4 & 2 \\
	2 & 2  & 5
\end{bmatrix}
=
-7.56
\begin{bmatrix}
	-0.98\\
	0.14\\
	0.13
\end{bmatrix}
\begin{bmatrix}
	-0.98\\
	0.14\\
	0.13
\end{bmatrix}
^{t}	
+2.44
\begin{bmatrix}
	-0.03\\
	-0.77\\
	0.63
\end{bmatrix}
\begin{bmatrix}
	-0.03\\
	-0.77\\
	0.63
\end{bmatrix}
^{t}
+7.12
\begin{bmatrix}
	0.19\\
	0.61\\
	0.76
\end{bmatrix}
\begin{bmatrix}
	0.19\\
	0.61\\
	0.76
\end{bmatrix}
^{t}
\]
%}$

\[
=
\begin{bmatrix}
	-7.26 & 1.04 & 0.96 \\
	1.04 & -0.15 & -0.14 \\
	0.96 & -0.14  & -0.14
\end{bmatrix}
+
\begin{bmatrix}
	0.00 & 0.06 & -0.05 \\
	0.06 & 1.45 & -1.18 \\
	-0.05 & -1.18  & 0.97
\end{bmatrix}
+
\begin{bmatrix}
	0.26 & 0.83 & 1.03 \\
	0.83 & 2.65 & 3.30 \\
	1.03 & 3.30  & 4.11
\end{bmatrix}
\]

The eigenvalues sorted in the descending order of absolute values are given by $-7.56$, $7.12$, and $2.44$ respectively. The corresponding eigenvectors are given by $[-0.98,0.14,0.13]^{t}$,$[0.19,0.61,0.76]^{t}$, and $[-0.03,-0.77,0.63]^{t}$. Thus, we could rewrite the optimization problem $max \; x^{t}Qx$ as $max \; \sum_{i=1}^{n} \lambda_i (x^{t} c_i)^2 $ which reduces to $max \; \sum_{i=1}^{n} \lambda_i \sum_{j=1}^n (x^j c_i^j)^2 $. Note that the index $j$ runs over all variables for specific eigenvectors $c_i$. Also, subterms like $(x^{t} c_i)^2 $ involves outer product of a vector by itself and could lead to increase in the overall density of the resulting sum. 

The theoretical properties and applications of eigendecomposition such as image processing are well studied. According to \cite{smith2002tutorial}, eigenvectors can only be determined for square $nXn$ matrices with $n$ dimensions. All the eigenvectors of a matrix are orthogonal and normalized such that length is $1$ unit. More importantly, eigenvectors identify the principal components of a matrix. These components characterize the data, which could be compressed using some top components. However, information loss occurs in these cases of dimensionality reduction. The eigenvector with the highest eigenvalue is defined as the principal component of the dataset. It captures the most significant relationship amongst the data dimensions.

By sorting the eigenvalues in descending order of absolute values, components are sorted in order of significance. Dimensionality reduction could be accomplished by dropping components of lower significance. In this paper, we emphasize the first $k$ eigenvectors in this sorted list to guide the search towards areas of importance. Note that the sign of the eigenvalue is also important since multiplying a vector by $-1$ reverses the direction of the eigenvector. On the other hand, the magnitude of the eigenvalue could bias the solution landscape if it is too big. Thus, we limit the impact of the magnitude of the eigenvalue through a non-negative parameter $M$. To summarize, we augment the original $Q$ matrix with a penalty or reward term dependent on $M$, the sign of the eigenvalue, and the outer product of the eigenvector with itself.  A related eigenvalue decomposition based branch and bound algorithm utilizing semidefinite relaxation was found effective for Quadratically Constrained Quadratic Programming instances (\cite{lu2017eigenvalue,lu2019sensitive}). However, our approach employs a simple transformation of the $Q$ matrix. %Note that a related eigenvalue decomposition based branch and bound algorithm have has found effective for Quadratically Constrained Quadratic Programming instances (\cite{}).

Heuristic solvers like tabu search and path relinking benefit from these augmented terms by implicitly (not directly) guiding the search towards promising areas in the solution landscape. In the above example, by considering the top most eigenvalue (sorted in descending order of absolute values), the contribution term is given by $1*(-1)*(-0.98 x_1 + 0.14 x_2 + 0.13 x_3)^2$ assuming $M=1$. Note that the sign of the corresponding eigenvalue $-7.56$ is $-1$. This is equivalent to adding $1*(-1)*[-0.98,0.14,0.13]^{t}*[-0.98,0.14,0.13]$ to the original $Q$ matrix. Thus, the original $Q$ matrix is transformed into $Q^{\prime}$ as follows:
\[
Q^{\prime} =
\begin{bmatrix}
	-7 & 2 & 2 \\
	2 & 4 & 2 \\
	2 & 2  & 5
\end{bmatrix}
- 1
\begin{bmatrix}
	-0.98\\
	0.14\\
	0.13
\end{bmatrix}
\begin{bmatrix}
	-0.98\\
	0.14\\
	0.13
\end{bmatrix}
^{t}
=
\begin{bmatrix}
	-7 & 2 & 2 \\
	2 & 4 & 2 \\
	2 & 2  & 5
\end{bmatrix}
- 1
\begin{bmatrix}
	0.96 & -0.14 & -0.13 \\
	-0.14 & 0.02 & 0.02 \\
	-0.13 & 0.02  & 0.02
\end{bmatrix}
\]
\[
= 
\begin{bmatrix}
	-7.96 & 2.14 & 2.13 \\
	2.14 & 3.98 & 1.98 \\
	2.13 & 1.98  & 4.98
\end{bmatrix}
\]

\section{Model}

The base problem of interest is $max \; x^{t}Qx$. We propose a transformation of the base problem into $max \; x^{t}Q'x = max \; x^{t}Qx + \sum_{i=1}^{k} M*sign(\lambda_i) * (x^{t} c_i)^2$. The index $i$ is based on the sorted list of absolute eigenvalues in descending order. Note that the augmented term is quadratic and fits into the QUBO modeling framework. Thus, $max \; x^{t}Q'x$ could be easily handled by any QUBO solver. We will compare the performance of a tabu-search based QUBO solver and CPLEX on both the base problem and the transformed problem in Section 3. Note that our proposed technique classifies as a preprocessing approach for QUBO similar to \cite{verma2020optimal}.

The magnitude of the penalty coefficient $M$ could be varied to emphasize the degree of importance. However, a large $M$ biases the solution landscape heavily towards the direction of the eigenvector. On the other hand, small values of $M$ are inconsequential and might not lead to any difference in the search process. Similarly, we can decide the number of top eigenvectors ($k$) to be considered in the contribution term. The search process will be overwhelmed if we consider a lot of eigenvectors. Moreover, the eigenvectors with small eigenvalues could be safely ignored since they do not contain a lot of information. We will perform sensitivity analysis on both parameters $M$ and $k$ in Section 3.

\section{Computational Experiments and Results}

For computational experiments we utilize the ORLIB datasets (\cite{beasley1990or}) and Palubeckis datasets (\cite{palubeckis2004multistart}). The former contains datasets with 100, 250, 500 and  1000 nodes, and integral coefficients of the $Q$ matrix are randomly chosen from the uniform distribution [-100,100]. These datasets are sparse with linear and quadratic density of 10\%. The latter contains larger and denser (10-100\% dense) datasets with 3000, 4000, 5000, 6000 and 7000 nodes. The uniform distribution of the coefficients of the $Q$ matrix is expected to generate a distribution of eigenvalues that are also uniform, whereas a non random $Q$ distribution, e.g. clustered $Q$, is expected to have an eigenvalue distribution with noticeable peaks corresponding to the clusters. We will investigate computational impact of the transformation on both types of $Q$ matrices in this section.

The algorithms were implemented in Python 3.6. The experiments were performed on a 3.40 GHz Intel Core i7 processor with 16 GB RAM running 64 bit Windows 7 OS. Smaller datasets up to 1000 nodes were solved using CPLEX 12.6.3. while the bigger datasets utilize a path relinking and tabu search based QUBO solver based on \cite{lewis2021qfold}. In this heuristic, the effect of a single bit-flip is quickly evaluated in $O(1)$. This enables the selection of the best variable in a neighborhood with the greatest effect in $O(n)$. This process is complemented by path relinking, which is implemented as a greedy search of the restricted solution space defined by the difference bits. Variables selected to be flipped are assigned to a tabu list to prevent cycling. When stuck in a local optima, a random restart is executed. We will refer to this heuristic as PRlocal in the forthcoming discussion.

First, we will investigate the distribution of eigenvalues for one of the ORLIB instances with 1000 nodes, $bqp\_1000\_1$. The results are summarized in Figure 1. We see a truncated normal distribution with the peak around zero. This behavior is consistent with other ORLIB and Palubeckis datasets. Note that these components with smaller absolute eigenvalues could be removed from the $Q$ matrix according to the associated eigendecomposition in Section 1. However, we will focus on the components with top $k$ absolute eigenvalues according to Section 2. Herein, the biggest such value is around 1000.  %Graph of distribution of eigenvalues for sparse (bqp1000) vs dense (3000 and 1)

\begin{figure}[htbp]
	\centerline{\includegraphics[scale=0.65]{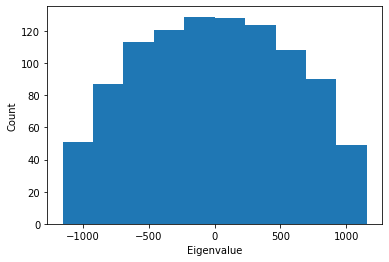}}
	\caption{Frequency distribution of eigenvalues for $bqp\_1000\_1$}
	\label{fig:fig01}
\end{figure}

Next, we will discuss the performance of the transformed model on CPLEX on smaller datasets containing fewer than 1000 nodes. We allow 100 seconds for each experimental run. Further, we utilized default CPLEX settings except \textit{qtolin}. We disabled the linearization capability of CPLEX by setting the \textit{qtolin} parameter to zero after noticing all the allocated time was invested in generating linear constraints for the binary quadratic problem. Note that we use $k=1$ as baseline parameter values since $k=1$ identifies the principal component having the most significant relationship between data dimensions (\cite{smith2002tutorial}). Moreover, the leading eigenvector is used prominently in various research papers (\cite{sarkar1998quantitative,perona1998factorization}). Our decision on $M=100$ was driven by the fact that eigenvectors are normalized (the individual entries are less than 1). Hence, the outer product of the eigenvector by itself would lead to off-diagonal entries, which are the multiple of two numbers less than 1. Multiplying by $100$ creates a contribution term which is significant but not very large relative to the entries in the $Q$ matrix. In other words, we want to perturb the original $Q$ matrix by a small amount such that the search process is not overwhelmed by the extraneous terms.

For the smaller problems with 100, 250 and 500 nodes, there was no observed difference in the best solution obtained by CPLEX for both the base problem and transformed problem described in Section 2. These datasets are smaller which allows CPLEX to quickly solve these problems within 100 seconds. For the 1000 node instances, we report the average percentage effect of the transfomation compared to the base problem across all datasets in the Improvement column of Table 1. In general, as $M$ and $k$ increases, the performance of the transformation deteriorates. %We observe that the transformation improves the solution in 6 out of 10 instances. The average percentage improvement across all datasets is 0.069\%.

Next, we analyze the impact of changing the parameters $M$ and $k$ on the computational performance. For this purpose, we experimented with $k \in \{1,2,5,10,20,25\}$ and $M \in \{100,200,300,400,500\}$ in Table 1. We notice that $M=100$ and $200$ lead to similar performance profiles. However, as the value of $M$ increases, the performance deteriorates and the transformation leads to reduction in solution quality. Utilizing up to top $10$ eigenvectors in the transformation leads to a gain of around 0.02\%. However, this gain diminishes as we consider more eigenvectors in the transformation.

\begin{table}
	\centering
	\caption{Impact of transformation on 1000 node instances  (CPLEX, $M=100,k=1$)}
	\begin{tabular}{l|lllll}
		\hline
		&  \multicolumn{5}{c}{Average Improvement \% Over Base Problem}  \\ \hline
		k & M=100 & M=200 & M=300 & M=400 & M=500 \\ \hline
		1 & 0.069 & 0.031 & -0.015 & -0.169 & -0.362 \\ 
		2 & 0.076 & 0.055 & -0.031 & -0.236 & -0.955 \\ 
		5 & 0.053 & -0.076 & -0.26 & -0.958 & -2.431 \\ 
		10 & 0.064 & -0.079 & -0.761 & -2.452 & -3.889 \\ 
		20 & 0.056 & -0.426 & -1.431 & -3.017 & -4.343 \\ 
		25 & -0.022 & -0.555 & -1.753 & -3.115 & -4.658 \\
	\end{tabular}
\end{table}

Similar to \cite{verma2020penalty}, we analyze the performance of the solver on the 1000 node instances through autocorrelation coefficient $\xi$. The higher value of $\xi$ is related to a lower number of local optima, which leads to better performance of any local search solver. We observe that the $\xi$ value for smaller values of $k$ and $M$ is somewhat higher. Thus, the gain in performance could be explained using a statistic from the search landscape analysis.

\begin{table}
	\centering
	\caption{Impact on solution landscape for 1000 node instances}
	\begin{tabular}{l|lllll}
		\hline
		& \multicolumn{5}{c}{Autocorrelation Coefficient ($\xi$)}  \\ 
		\hline
		k & M=100 & M=200 & M=300 & M=400 & M=500 \\ 
		\hline
		1 & 356.1 & 356.2 & 354.4 & 355.5 & 354.4 \\ 
		2 & 354.4 & 353.0 & 352.8 & 353.4 & 352.8 \\ 
		5 & 354.3 & 352.7 & 352.5 & 352.6 & 352.0 \\ 
		10 & 352.9 & 352.4 & 352.4 & 352.0 & 349.6 \\ 
		20 & 352.6 & 352.0 & 350.8 & 349.2 & 348.0 \\ 
		25 & 351.1 & 351.8 & 348.5 & 346.7 & 346.4 \\ 
	\end{tabular}
\end{table}

For the bigger Palubeckis instances, we performed sensitivity analysis on $M$ and $k$ using the baseline values of $M=100$ and $k=1$, and the results are presented in Figure 2. We observed similar performance profiles based on Average Improvement \% across other bigger instances. Thus, $k \in \{1,2,5,10\}$ and $M \in \{100,200\}$ are also acceptable choices for bigger datasets.

\begin{figure}%
	\centering
	\subfloat[\centering Sensitivity analysis of $k$]{{\includegraphics[width=6.25cm]{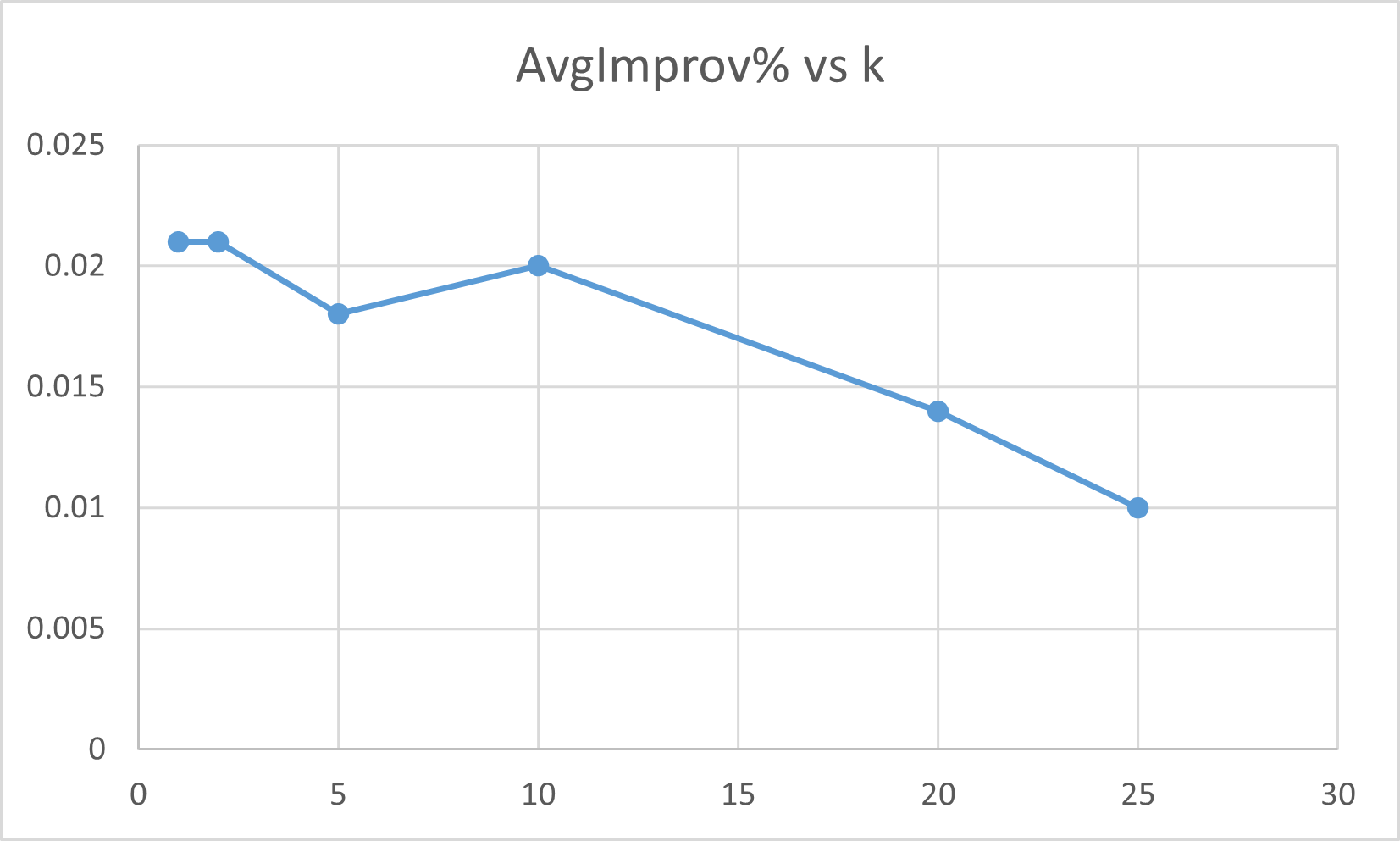} }}%
	\qquad
	\subfloat[\centering Sensitivity analysis of $M$]{{\includegraphics[width=6.25cm]{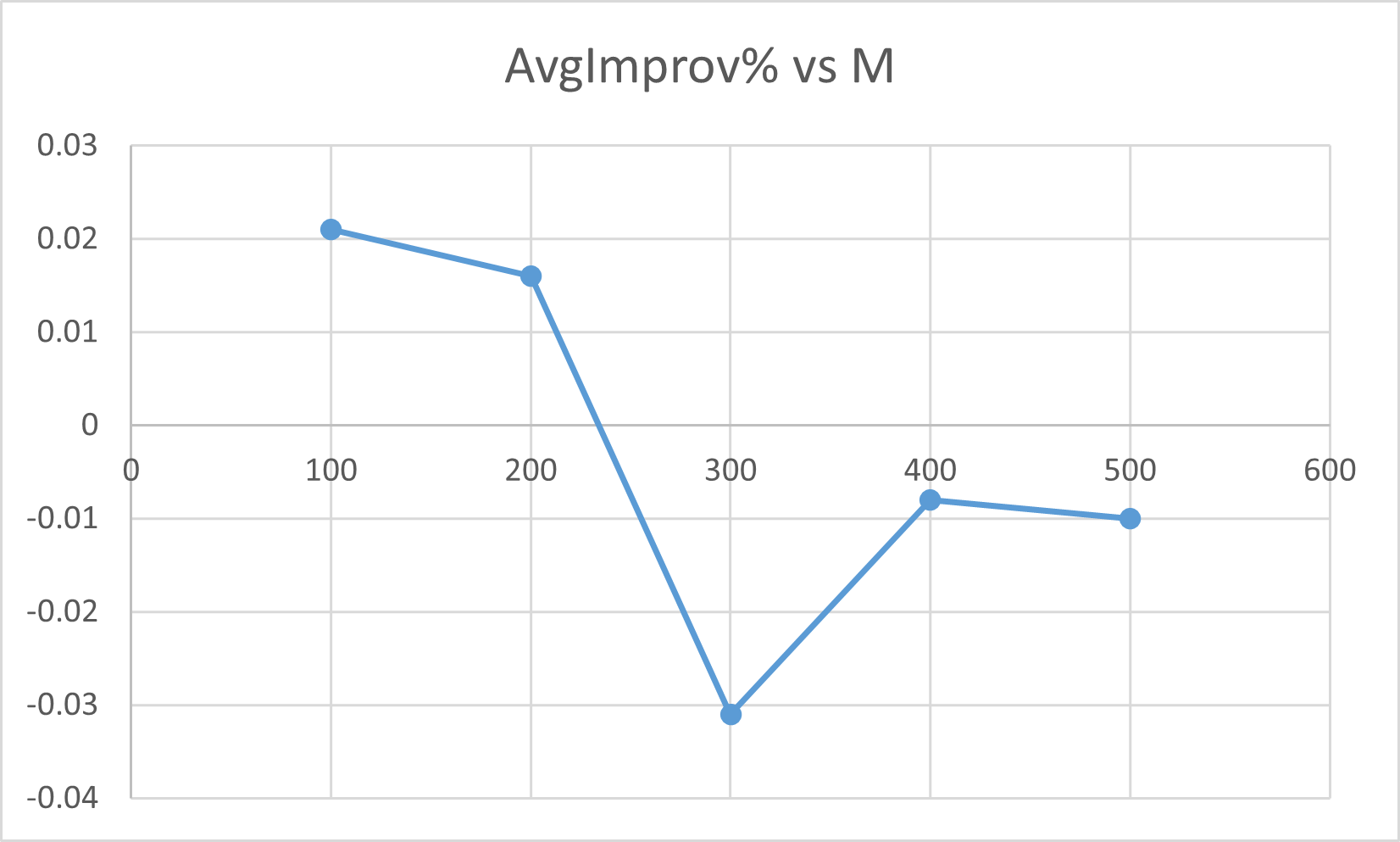} }}%
	\caption{Results for bigger Palubeckis instances}%
	\label{fig:example}%
\end{figure}

The results on Average Improvement \% across all bigger instances using $M=100$ and $k=1$ are presented in Table 2. The transformation improves or matches the base problem in 27 out of 31 instances. We want to mention that the additional contribution term due to eigenvalues i.e. $M*sign(\lambda_i) * (x^{t} c_i)^2$ leads to the creation of a dense matrix because of the inner product. Thus, the transformation is more effective on Palubeckis instances which are already dense and the adjustment does not create new off-diagonal entries.

\begin{table}
	\centering
	\caption{Impact of transformation on bigger instances (PRlocal, $M=100,k=1$)}
	\begin{tabular}{l|l|l|l}
		\hline
		Instance Type &                Nodes & Avg Improvement & Avg Improvement \% \\ \hline
		ORLIB &                         2500 & 2478 &                 0.021 \\ \hline
		Palubeckis &                 3000 & 385 &                         0.017 \\ 
		&                                 4000 & 16 &                         0.011 \\ 
		&                                 5000 & 5367 &                 0.049 \\ 
		&                                 6000 & 5018 &                 0.065 \\ 
		&                                 7000 & 3875 &                 0.031 \\ \hline
	\end{tabular}
\end{table}

Based on our expectation that $Q$ with a few dominant eigenvalues will be more amenable to this technique, we investigate the impact of the transformation on such matrices. For this purpose, we study the Maximum Diversity Problem (MDP) instances discussed in \cite{marti2013heuristics}. The optimization problem of interest $max \sum_{i,j,i \neq j} d_{ij} x_i x_j \; s.t. \; \sum_i x_i = m, \; x_i \in \{0,1\}, \; i \in [1,n]$ comprises of selecting a subset of $m$ elements from $n$ variables such that the sum of distance (or equivalent diversity measure) between chosen elements is maximized. The QUBO transformation incorporates the quadratic penalty term given by $-P(\sum_i x_i - m)^2$ resulting in a $Q$ matrix where $P$ represents a non-negative penalty coefficient. The details of the QUBO transformation for MDP are discussed in \cite{ghosh1996computational}. The set of instances identified as MDPLIB consists of 315 instances belonging to three categories SOM, GKD and MDG. These categories have small, medium and large size instances respectively. 

The eigenvalue distribution for $SOM\_a\_21\_n100\_m10$ is detailed next. We utilize $P=10, m=10$ and $n=100$ for the resulting QUBO model. Figure 3 clearly indicates that there is a dominant eigenvalue around -500 while the remaining eigenvalues are concentrated around 200. Other MDPLIB instances follow a similar trend. Thus, the principal component given by the largest absolute eigenvalue clearly dominates other components. Hence, guiding the search towards this principle component might have an impact on the computational performance of the solver.

\begin{figure}[htbp]
	\centerline{\includegraphics[scale=0.65]{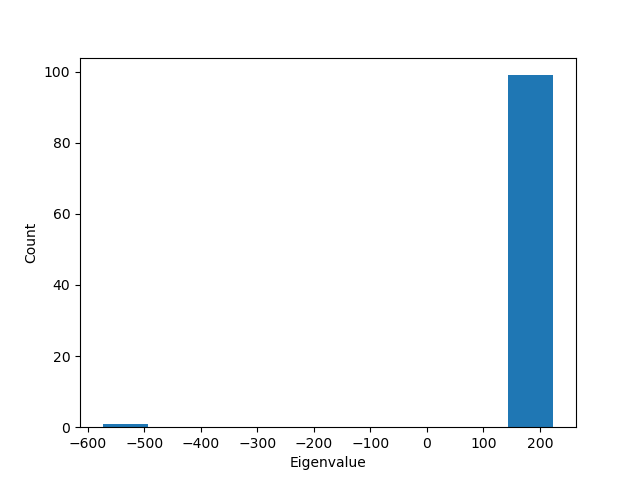}}
	\caption{Frequency distribution of eigenvalues for $SOM\_a\_21\_n100\_m10$}
	\label{fig:fig02}
\end{figure}

Again, we utilize $M=100$, $k=1$ and $t=100$ seconds for these experiments. Similar to smaller ORLIB QUBO benchmark instances, SOM and GKD have atmost 500 variables. The eigenvalue decomposition has no improvement on solution quality obtained by PRlocal because of smaller size. On the other hand, larger MDG instances see a significant increase in solution quality as a result of the transformation. Note that the Average Improvement \% in Table 4 increases as the size of the QUBO instance increases. Overall, the eigendecomposition improves or matches the best solution found by the original problem in all but 17 instances. Thus, our technique is found to more effective in $Q$ matrices with dominant eigenvalues.

\begin{table}
	\centering
	\caption{Impact of transformation on MDPLIB MDG instances (PRlocal, $M=100,k=1$)}
	\begin{tabular}{l|l|l}
		\hline
		Instance Type &                Nodes &  Avg Improvement \% Over Base Problem\\ \hline
		MDG-a &                         500,2000 &                 1.0 \\ 
		MDG-b &                 500,2000 &                          6.0 \\ 
		MDG-c &                                 3000 &                          9.7 \\  \hline
	\end{tabular}
\end{table} 

%CPLEX

%PRlocal

\section{Conclusions and Future Research}

We present a learning-based technique that modifies the $Q$ matrix based on the principal components of the $Q$ matrix obtained from eigenvalue decomposition. Testing indicates this technique leads to improvement in solution quality for benchmark QUBO instances especially when there are a few dominant eigenvalues. Future work involves adjusting $M$ based on the eigenvalue for each component in the decomposition and dimensionality reduction of the $Q$ matrix.

%\section{Conclusions} \label{conc}
\section{Data Availability Statement}
The datasets generated during and/or analyzed during the current study are available from the corresponding author on reasonable request.

\section{Conflict of Interest Statement}
On behalf of all authors, the corresponding author states that there is no conflict of interest.

\bibliography{decompose}
%{\footnotesize%\scriptsize%\small
%\bibliography{rFilp}}


%{\footnotesize \bibliography{qsubs} }  % name your BibTeX data base
%% \bibliographystyle{elsarticle-harv} 

%\begin{thebibliography}{00}


%\bibitem{}

%\end{thebibliography}
\end{document}